\documentclass{amsart}

\usepackage{amsmath}
\usepackage{verbatim}
\usepackage{upref}
\usepackage{amsfonts,amssymb}
\usepackage{youngtab}

\newtheorem{theorem}{Theorem}[section]
\newtheorem{lemma}[theorem]{Lemma}
\newtheorem{proposition}[theorem]{Proposition}
\newtheorem{corollary}[theorem]{Corollary}

\newtheorem{main lemma}[theorem]{Main Lemma}
\newtheorem{conjecture}[theorem]{Conjecture}

\theoremstyle{definition}
\newtheorem{definition}[theorem]{Definition}

\newtheorem{remark}[theorem]{Remark}

\newcommand{\N}{\mathbb{N}}

\newcommand{\Z}{\mathbb{Z}}
\newcommand{\FF}{F\langle X|G\rangle}

\newcommand{\gk}{\text{\rm GKdim}}

\newcommand{\ep}{\text{\rm exp}}
\newcommand{\s}{\mathcal{S}}

\begin{document}\title[Graded GK dimension of the Grassmann algebra]{$\Z_2$-graded Gelfand-Kirillov dimension of the Grassmann algebra}\author{Lucio Centrone}\address{IMECC, Universidade Estadual de Campinas, Campinas (SP) Brazil, Rua Sergio Buarque de Holanda 651,13083-859}\email{centrone@ime.unicamp.br}\thanks{Partially supported by FAPESP 2013/06752-4 }\begin{abstract}We consider the infinite dimensional Grassmann algebra $E$ over a field $F$ of characteristic 0 or $p$, where $p>2$, and we compute its $\Z_2$-graded Gelfand-Kirillov dimension as a $\Z_2$-graded PI-algebra.\end{abstract}\maketitle

\section{Introduction}
The infinite dimensional Grassmann algebra $E$ plays an important role in the theory of algebras with polynomial identities (PI-algebras). Not to mention that the Grassmann algebra fits in the description of polynomial identities of any PI-algebra over a field of characteristic 0. It is well known that each polynomial identity of the Grassmann algebra $E$ is a consequence of the triple commutator $[x_1,x_2,x_3]$ (Krakovsky and Regev \cite{krr1}). Moreover $E$ is a superalgebra endowed with its natural grading. If $A=A^0\oplus A^1$ is a superalgebra, we denote by $G(A)$ the Grassmann envelope of $A$ that is $G(A)=(E^0\otimes A^0)\oplus (E^1\otimes A^1)$. If $A$ is a PI-algebra, then the ideal of polynomial identities of $A$, i.e. the $T$-ideal of $A$ or $T(A)$, is generated by the Grassmann envelope of a finite dimensional algebra (see the monograph of Kemer \cite{kem1}). 

We recall that if $X$ is an infinite countable set of indeterminates and $F$ is a field, a $T$-ideal of $F\langle X\rangle$, the free-associative algebra freely generated by $X$, is an ideal invariant under all the endomorphisms of $F\langle X\rangle$. It is well known that $T(E)$ is a $T$-prime ideal, i.e., an ideal that is prime as a $T$-ideal. This is of crucial importance in light of the structure theory of $T$-ideals developed by Kemer. Of course the $T$-ideal of the Grassmann algebra is not the only $T$-prime ideal. In fact the non-trivial $T$-prime ideals are: $T(M_n(F))$, the $T$-ideal of the algebra of $n\times n$ matrices with entries from a field $F$; $T(M_n(E))$, the $T$-ideal of the algebra of $n\times n$ matrices with entries from $E$ and $T(M_{a,b}(E))$, the $T$-ideal of the block-matrix algebra with entries from $E^0$ in the upper left $a\times a$ corner and in the lower right $b\times b$ corner and with entries from $E^1$ elsewhere.

The deep structure theory of $T$-ideals is no more valid in positive characteristic. Nevertheless the interest of mathematicians in the polynomial identities of $E$ in positive characteristic is proved in a number of papers. Let $F$ be the ground field of $E$, then in the case $F$ is finite, Regev showed that the triple commutator is not enough to describe the whole $T$-ideal of $E$ (see \cite{reg2}). By the way
when the field is infinite the situation is close to the case of characteristic 0. In fact in \cite{gik1} Giambruno and Koshlukov show that all the polynomial identities of $E$ follow from the triple commutator. In light of the last result the authors gave a partial answer to a question posed by Kemer in \cite{kem2}, i.e., does $E$ satisfy all the polynomial identities of $M_{\frac{p+1}{2}}(F)$, where $p$ is the characteristic of $F$? In particular they proved that in case $p=3$ $E$ does not satisfies the identities of $M_2(F)$ unless the finite dimensional Grassmann algebras of dimension 4 and 5.

In the last five years the interest moved to the graded case and, in particular, when $E$ is endowed with a $\Z_2$-grading. Both in the case $F$ is of characteristic 0 and positive, we have a complete description of the $\Z_2$-graded polynomial identities of $E$. In \cite{did1} Di Vincenzo and da Silva described all the possible $\Z_2$-homogeneous gradings over $E$. It turns out there are three different classes of $\Z_2$-gradings for $E$. In order to differentiate the three gradings they use a different notation to indicate $E$ in each of the cases, i.e., $E_{k^*}$, $E_\infty$, $E_k$. They also gave a complete list of identities which generate the set of $\Z_2$-graded identities of $E_d$  ($d\in\{k^*,\infty,k\}$). Using similar techniques the author gave a complete list of identities which generate the set of $\Z_2$-graded identities of $E_d$ in case of positive characteristic (see \cite{cen3}).

Following this line of research we want to focus on the $\Z_2$-graded relatively-free algebra $A:=F\langle X\rangle/I$, where $I$ is the ideal of $\Z_2$-graded identities of $E_d$ and $X=\{y_1,\ldots,y_m,z_1,\ldots,z_m\}$ is a set of homogeneous variables such that the $y_i$'s belong to the 0 component and the $z_j$'s to the component of degree 1. In particular we exhibit a multihomogeneous basis for $A$. We use such a result in order to compute the Gelfand-Kirillov (GK) dimension of $A$, i.e. the $\Z_2$-graded GK dimension of $E_d$, and we write down explicitely the Hilbert series of $A$. The three main theorems deal with the GK dimension of $A$ for $E_d$. As a consequence we have the following result:

\begin{theorem}
Let $F$ be a field of positive characteristic. Then the GK dimension of $A$ is $m$ except in the case $d=\infty$.
\end{theorem}

\section{Basic tools}

All algebras we refer to are associative and unitary unless explicitely written. Moreover all the field are to be considered infinite and of characteristic different from 2.

\begin{definition}
Let $(G,\cdot)=\{g_1,\ldots,g_r\}$ be any group and let $F$ be a field. If $A$ is an $F$-algebra, we say that $A$ is a
$G$-graded algebra if there are subspaces $A^g$ for each $g\in G$ such that \[A=\bigoplus_{g\in G}A^g \ \textrm{and} \
A^gA^h\subseteq A^{gh}.\] If $0\neq a\in A^g$ we say that $a$ is \textit{homogeneous of $G$-degree $g$} or simply that $a$ has degree $g$,
and we write $\deg(a)=g$.\end{definition}

One defines \textit{$G$-graded subspaces} of $A$, \textit{$G$-graded $A$-modules}, \textit{$G$-graded homomorphisms} and so on, in a standard way (see for
example \cite{dre1} or \cite{ddn1} for more details).

Let $\{X^{g}\mid g \in G\}$ be a family of disjoint countable sets. Put $X=\bigcup_{g\in G}X^{g}$ and denote by $\FF$ the free-associative algebra freely generated by the set $X$. An indeterminate $x\in X$ is said to be of \textit{homogeneous $G$-degree $g$}, written $\deg(x)=g$,
if $x\in X^{g}$. We always write $x^{g}$ if $x\in X^{g}$. The homogeneous $G$-degree of a monomial $m=x_{i_1}x_{i_2}\cdots x_{i_k}$ is defined to be
$\deg(m)=\deg(x_{i_1})\cdot\deg(x_{i_2})\cdot\cdots\cdot\deg(x_{i_k})$. For every $g \in G$, we denote by $\FF^g$ the subspace of $\FF$ spanned by all the
monomials having homogeneous $G$-degree $g$. Notice that $\FF^g\FF^{g'}\subseteq \FF^{gg'}$ for all $g,g' \in G$. Thus \[\FF=\bigoplus_{g\in G}\FF^g\] is a $G$-graded algebra. The elements of the $G$-graded algebra $\FF$ are referred to as \textit{$G$-graded polynomials} or, simply, \textit{graded polynomials}. We refer to the multihomogeneous degree of a polynomial as the usual multidegree of $\FF$ when naturally graded by $\Z$.

\begin{definition} If $A$ is a $G$-graded algebra, we denote by $T_{G}(A)$ the intersection of the kernels of all
$G$-graded homomorphisms $\FF\rightarrow A$. Then $T_{G}(A)$ is a graded two-sided ideal of $\FF$ and its elements are called \textit{$G$-graded polynomial identities} of the algebra $A$.\end{definition}
Notice that $T_{G}(A)$ is stable under the action of any $G$-graded endomorphism of the algebra $\FF$.
Any $G$-graded ideal of $\FF$ which verifies such property is said to be a \textit{$T_{G}$-ideal}. Clearly,
any $T_{G}$-ideal $I$ is the ideal of the $G$-graded polynomial identities of the graded algebra
$\FF/I$. Note also that for a $G$-graded algebra $A$, the quotient algebra $\FF/T_G(A)$ is the
relatively-free algebra for the variety of $G$-graded algebras generated by $A$.
Let $A$ be any $G$-graded PI-algebra over $F$, where $G$ is any finite group of order $s$. We shall denote by $F_k^{G}(A)$ the relatively-free algebra \[{F\langle x_1^{g_1},\ldots,x_k^{g_1},\ldots,x_1^{g_s},\ldots,x_k^{g_s}\rangle}/({F\langle x_1^{g_1},\ldots,x_k^{g_1},\ldots,x_1^{g_s},\ldots,x_k^{g_s}\rangle\cap T_{G}(A)}).\] The latter is called \textit{relatively-free $G$-graded algebra} of $A$ in $k$ variables. In the case $G$ is the trivial group, the previous definition coincides with the definition of relatively-free algebra for PI-algebras. Moreover, we shall use the symbol $f\equiv g$ in order to say that the graded polynomials $f$ and $g$ are equal modulo the graded identities of a certain graded algebra.

In this paper we are going to deal with $\Z_2$-gradings over the infinite dimensional Grassmann algebra. Sometimes we use the word \textit{superalgebra} instead of $\Z_2$-graded algebra.

\begin{definition}
Let $V=\{v_1,v_2,\cdots\}$ be an infinite countable set, then we denote by $E=E(V)$ the Grassmann algebra generated by $V$, i.e., $F\langle V\rangle/I$, where for each $i$, $e_i=v_i+I$ and $I$ is the ideal generated by $\{v_iv_j+v_jv_i|i,j\in\N\}$. It is well known that $B_E=\left\{e_{i_1}e_{i_2}\cdots e_{i_n}\mid n\in \N, {i_1}<{i_2}<\cdots<{i_n}\right\}$ is a basis of $E$ as a vector space over $F$. Moreover we say $e_{i_1}\cdots e_{i_l}$ is a basis element of $E$ of \textit{length} $l$ with \textit{support} $\{e_{i_1},\ldots,e_{i_l}\}$. 
\end{definition}

Let us consider the map $\varphi:V\rightarrow \Z_2$ such that $x_i\mapsto 1$. The map $\varphi$ gives out a $\Z_2$-grading over $E$ called \textit{natural grading}. In this case, let $E^0$ be the homogeneous component of $\Z_2$-degree 0 and let $E^1$ be the component of degree 1. It is easy to see that $E^0$ is the center of $E$ and $ab + ba = 0$ for all $a, b\in E^1$. This means that $E$ satisfies the following graded polynomial identities:
$[y_1, y_2]$, $[y_1, z_1]$, $z_1z_2 + z_2z_1$.
Now, let us consider the $\Z_2$-gradings over $E$ induced by the maps  $\deg_{k*}$, $\deg_\infty,$ and $\deg_k,$ defined respectively by

$$\deg_{k*}(e_i)=\left\{\begin{array}{ll}
\text{\rm 1 for $i = 1,\ldots, k$}\\
\text{\rm 0 otherwise},\end{array}\right.$$

$$\deg_{\infty}(e_i)=\left\{\begin{array}{ll}
\text{\rm 1 for $i$ odd}\\
\text{\rm 0 otherwise},\end{array}\right.$$

$$\deg_k(e_i)=\left\{\begin{array}{ll}
\text{\rm 0 for $i = 1,\ldots, k$}\\
\text{\rm 1 otherwise}.\end{array}\right.$$

From now on we shall denote by $E_{k^*}$, $E_\infty$, $E_k$ the Grassmann algebra endowed with the $\Z_2$-grading induced by the maps $\deg_{k*}$, $\deg_\infty,$ and $\deg_k$. We denote by $E_d$ any of the superalgebras $E_{k^*}$, $E_\infty$, $E_k$ without distinguish them.

\begin{remark}\label{possiblechoose}
We note that we may always choose elements of $E_d^0$ of the form $e_{i_1}e_{i_2}+\cdots+e_{i_{2a-1}}e_{i_{2a}}$ for each $a\in\N$.
\end{remark}

We introduce a special type of polynomials that turned out to be crucial when describing the $\Z_2$-graded identities of $E$ (see \cite{did1}).

\

Let $f=z_{i_1}^{r_{i_1}}\cdots z_{i_s}^{r_{i_s}}[z_{j_1},z_{j_2}]\cdots [z_{j_{t-1}},z_{j_t}]$ and consider the set \[\s=\{\text{different homogeneous variables appearing in $f$}\}\subseteq\{z_1,\ldots,z_m\}.\] If $h=|\s|,$ then $\s=\{z_{i_1},\ldots,z_{i_h}\}.$  Notice also that $f$ is linear in the commutators.

We consider now \[T = \{j_1,\ldots, j_t\}\subseteq\s\] and let us denote the previous polynomial by \[f_T (z_{i_1},\cdots, z_{i_h}).\] 

\begin{definition} For $m\geq2$ let $$g_m(z_{i_1},\ldots, z_{i_h}) =\sum_{\begin{array}{cc}
T\\
\text{\rm $|T|$ even}\end{array}}(-2)^{-\frac{|T|}{2}}f_T(z_{i_1},\ldots, z_{i_h}),$$
moreover put
$g_1(z_1) = z_1.$\end{definition}

We have the following theorems due to the author \cite{cen3}. 

\begin{theorem}\label{t1} Let $k\in\N$ and $p>2$ be a prime. Let us consider the Grassmann algebra over a field of characteristic. If $p>k$, then all $\Z_2$-graded polynomial identities of $E_{k^*}$ are consequences of the graded identities: \[[x_1,x_2,x_3],\ z_1\cdots z_{k+1}.\] 
On the other side, if $p\leq k,$ all $\Z_2$-graded polynomial identities of $E_{k^*}$ are consequences of the graded identities: \[[x_1,x_2,x_3],\ z_1\cdots z_{k+1},\ z^p.\]
\end{theorem}

\begin{theorem}\label{t2} Let $p>2$ be a prime. Let us consider the Grassmann algebra over a field of characteristic. Then all $\Z_2$-graded polynomial identities of $E_\infty$ are consequences of the graded identities: \[[x_1,x_2,x_3],\ z^p.\]\end{theorem}

\begin{theorem}\label{t3}  Let $k\in\N$ and $p>2$ be a prime. Let us consider the Grassmann algebra over a field of characteristic. If $p>k$, then all $\Z_2$-graded polynomial identities of $E_k$ are consequences of the graded identities:\begin{itemize}
\item $[x_1,x_2,x_3],$
\item \text{\rm $[y_1,y_2]\cdots[y_{k-1},y_k][y_{k+1},x]$ (if $k$ is even)}
\item \text{\rm $[y_1,y_2]\cdots[y_{k},y_{k+1}]$ (if $k$ is odd)}
\item \text{\rm $g_{k-l+2}(z_1,\ldots,z_{k-l+2})[y_1,y_2]\cdots[y_{l-1},y_l]$ (if $l\leq k$, $l$ is even)}
\item \text{\rm $[g_{k-l+2}(z_1,\ldots,z_{k-l+2}),y_1][y_2,y_3]\cdots[y_{l-1},y_l]$ (if $l\leq k,$ $l$ is odd)}
\item \text{\rm $g_{k-l+2}(z_1,\ldots,z_{k-l+2})[z,y_1][y_2,y_3]\cdots[y_{l-1},y_l]$ (if $l\leq k,$ $l$ is odd).}
\end{itemize}
If $p\leq k$ we have to add to the list above the identity
\begin{itemize}
\item $z^p.$
\end{itemize}
\end{theorem}

We recall that the analog of the previous theorems in the case of characteristic 0 has been completely solved by Di Vincenzo and da Silva in \cite{did1} and most of the techniques used in \cite{cen3} are due to this work.

We denote by $E^*$ the Grassmann algebra without 1. Then we have the following result (see \cite{gik1} for more details).

\begin{proposition}\label{regev}
Let $F$ be a field of characteristic $p>0$, then $x^p$ is a polynomial identity for $E^*$.
\end{proposition}

By Proposition \ref{regev}, we have the following lemma.

\begin{lemma}\label{regevcons1}
Let $F$ be a field of characteristic $p>0$, $f\in F$ and $a\in E^*$ be a basis element of even length, then $(f+a)^p=f^p$.
\end{lemma}
\proof
Due to the fact that $a$ is in the center of $E$ and $F$ is of characteristic $p>0$, we have $(f+a)^p=f^p+a^p$, then the lemma follows by Proposition \ref{regev}.
\endproof

\begin{corollary}\label{regevcons2}
Let $F$ be a field of characteristic $p>0$ and $\alpha,\beta\in\N$ such that $\text{\rm $\alpha\equiv \beta($mod $p)$}$. If $a,b\in E^*$ are basis elements of even length, then $(f+a)^\alpha=f^{pk}(f+b)^\beta$ for some $k\in\Z$.
\end{corollary}

We have the following result.

\begin{lemma}\label{regevcons3}
Let $F$ be an infinite field of characteristic $p>2$ and $\alpha_1<\alpha_2<\cdots<\alpha_l\in\N$, then in the $\Z_2$-graded relatively-free algebra $F_m^{\Z_2}(E_d)$ the monomials $y^{\alpha_1},\ldots,y^{\alpha_l}$ are linearly independent. 
\end{lemma}

\section{The graded Gelfand-Kirillov dimension}

We recall briefly some definitions and results about the graded Gelfand-Kirillov dimension of a graded PI-algebra. We refer to the surveys \cite{dre2} of Drensky and \cite{cen4} of the author for more details.

\begin{definition} \label{def2} Let $G=\{g_1,\ldots,g_s\}$ be a finite group and $A$ be a $G$-graded PI-algebra. The \textit{$G$-graded Gelfand-Kirillov (GK) dimension} of $A$ in $k$ variables is \[{\rm GKdim}^G_k(A):={\rm GKdim}(F_k^G(A)),\] where the new $sk$ non-commutative variables are $x_1^{g_1},\ldots,x_1^{g_s},\ldots,x_k^{g_1},\ldots,x_k^{g_s}$. \end{definition}

Notice that we have defined the $G$-graded GK dimension of $A$ in $k$ variables as the GK dimension of $F_k^G(A)$ that is generated by $sk$ variables. 

\begin{definition}\label{defhilb} Let $G=\{g_1,\ldots,g_s\}$ be a finite abelian group of order $s$. Let $A$ be a $G$-graded PI-algebra and $F_{k_1,\ldots,k_s}$ be its relatively-free $G$-graded algebra in the variables $T_1:=\{{t_1}_1,\ldots,{t_1}_{k_1}\}$,$\ldots$,\ $T_s:=\{{t_s}_1,\ldots,{t_s}_{k_s}\}$, where the variables in $T_1$ are of $G$-degree $g_1$ and the variables in $T_s$ are of $G$-degree $g_s$. It is well known that $F_{k_1,\ldots,k_s}$ is a  $\Z^{k_1+\cdots+k_s}$-graded algebra. The formal power series \[H^{G}(A;T_1,\ldots,T_s)=\]\[\sum_n\dim F_{k_1,\ldots,k_s}^{(n_1,\ldots,n_{k_1},\ldots,n_{\sum_{i=1}^{{s-1}}k_i+1},\ldots,n_{\sum_{i=1}^{{s}}k_i})}{t_1}_1^{n_1}\cdots {t_1}_{k_1}^{n_{k_1}}\cdots {t_s}_1^{n_{\sum_{i=1}^{s-1}k_i+1}}\cdots {{t_{s}}_{k_s}}^{n_{\sum_{i=1}^{s}k_i}}\] is called \textit{$G$-graded Hilbert series} of $A$ in the sets of variables $T_1,\ldots,T_s$. \end{definition}

\begin{definition}\label{remarkable}The growth function of $F_{k_1,\ldots,k_s}$ with respect to the vector space \[\text{\rm $V=$span$_F\langle T_1\cup\cdots\cup T_s\rangle$}\] is  \[g_V(n)=\sum_{n_1+\cdots +n_{\sum_{i=1}^sk_i}=n}\dim_FF_{k_1,\ldots,k_s}^{(n_1,\ldots,n_{\sum_{i=1}^{{s}}k_i})}.\]
\end{definition}

In \cite{pro1} and in \cite{ber1}, we may find a complete description of the Gelfand-Kirillov dimension of the PI-algebras whose $T$-ideals are to be considered the structure-blocks of the $T$-ideals of PI-algebras, i.e., the \textit{verbally prime algebras}. 

\section{The graded GK dimension of $E_{k^*}$}
We consider the algebra $E_{k^*}$ and we show a basis of its $\Z_2$-graded relatively-free algebra. Then we compute its $\Z_2$-graded GK dimension.

 Let us consider $k\in\N$ and $p>2$ be a prime such that $p>k$ and $E_{k^*}$ as defined above over a field of characteristic 0 or $p$. By Theorem 10 of \cite{did1} and by Theorem \ref{t1} the superalgebra $E_{k^*}$ satisfies the graded identity $[x_1,x_2,x_3]$, then each polynomial in the relatively-free graded algebra of $E_{k^*}$ can be written as a linear combination of the following:
\[y_{i_1}^{\alpha_{i_1}}\cdots y_{i_l}^{\alpha_{i_l}}z_{j_1}^{\beta_{j_1}}\cdots z_{j_r}^{\beta_{j_r}}[z_{s_1},z_{s_2}]\cdots[z_{s_p},z_{s_{p+1}}][y_{t_1},y_{t_2}]\cdots[y_{t_q},y_{t_{q+1}}]\]
\[y_{i_1}^{\alpha_{i_1}}\cdots y_{i_l}^{\alpha_{i_l}}z_{j_1}^{\beta_{j_1}}\cdots z_{j_r}^{\beta_{j_r}}[z_{s_1},z_{s_2}]\cdots[z_{s_p},z_{s_{p+1}}][z_{s_{p+2}},y_{t_1}][y_{t_2},y_{t_3}]\cdots[y_{t_q},y_{t_{q+1}}],\] where $i_1<i_2<\cdots<i_l$, $j_1<j_2<\cdots<j_r$, $s_1<s_2<\cdots<s_{p+1}$ and $t_1<t_2<\cdots<t_{q+1}$. Due to the identity $z_1\cdots z_{k+1}$ we may suppose the total degree, with respect to the $z$'s, of the previous polynomials to be less than or equal to $k$. In the case $p\leq k$ we may suppose each $\beta_{j_l}$ to be strictly less than $p$. Let us call \textit{$S_1$-type polynomials} those polynomials of the first of the two previous cases. Then we shall call \textit{$S_2$-type polynomials} those polynomials which are not in $S_1$.  Let $P=y_{i_1}^{\alpha_{i_1}}\cdots y_{i_l}^{\alpha_{i_l}}z_{j_1}^{\beta_{j_1}}\cdots z_{j_r}^{\beta_{j_r}}[z_{s_1},z_{s_2}]\cdots[z_{s_p},z_{s_{p+1}}][y_{t_1},y_{t_2}]\cdots[y_{t_q},y_{t_{q+1}}]$. From now on we shall refer to \[y_{i_1}^{\alpha_{i_1}}\cdots y_{i_l}^{\alpha_{i_l}}z_{j_1}^{\beta_{j_1}}\cdots z_{j_r}^{\beta_{j_r}}\] as the \textit{polynomial part} of $P$ and to \[[z_{s_1},z_{s_2}]\cdots[z_{s_p},z_{s_{p+1}}][y_{t_1},y_{t_2}]\cdots[y_{t_q},y_{t_{q+1}}]\] as the \textit{commutator part} of $P$. 

Now we have the following easy lemma.

\begin{lemma}\label{useful1}Let $p>2$ be a prime, $r\in\N$ and let $F$ be a field of characteristic 0 or $p>r$. If $s=e_1+e_{2}e_3+\cdots e_{2l}e_{2l+1}$, then $s^r=0$ if and only if $r\geq l+2$.
\end{lemma}
\proof
We observe that $2e_1e_2e_3$ is a monomial of $s^2$, then $6e_1e_2e_3e_4e_5$ is a monomial of $s^3$ and so on. We have $s^r=r!e_1e_2e_3\cdots e_{2l}e_{2l+1}\neq0$ and $s^rs=0$. 
\endproof

\begin{lemma}\label{independence} The polynomials of the following type: 

\[y_{i_1}^{\alpha_{i_1}}\cdots y_{i_l}^{\alpha_{i_l}}z_{j_1}^{\beta_{j_1}}\cdots z_{j_r}^{\beta_{j_r}}[z_{s_1},z_{s_2}]\cdots[z_{s_p},z_{s_{p+1}}][y_{t_1},y_{t_2}]\cdots[y_{t_q},y_{t_{q+1}}]\]
\[y_{i_1}^{\alpha_{i_1}}\cdots y_{i_l}^{\alpha_{i_l}}z_{j_1}^{\beta_{j_1}}\cdots z_{j_r}^{\beta_{j_r}}[z_{s_1},z_{s_2}]\cdots[z_{s_p},z_{s_{p+1}}][z_{s_{p+2}},y_{t_1}][y_{t_2},y_{t_3}]\cdots[y_{t_q},y_{t_{q+1}}].\]
 are linearly independent modulo $T_{\Z_2}(E_{k^*})$.
\end{lemma}
\proof
Let $F$ be of characteristic 0. We will argue only for the $S_1$-type polynomials because in the other case the proof is analogous. Let us consider the set $S_1(n)$ of polynomials of $S_1$ \[y_{i_1}^{\alpha_{i_1}}\cdots y_{i_l}^{\alpha_{i_l}}z_{j_1}^{\beta_{j_1}}\cdots z_{j_r}^{\beta_{j_r}}[z_{s_1},z_{s_2}]\cdots[z_{s_p},z_{s_{p+1}}][y_{t_1},y_{t_2}]\cdots[y_{t_q},y_{t_{q+1}}]\] of total degree $n$. Suppose we have a linear combination $\sum_P\gamma(P)P\equiv0$ of elements of $S_1(n)$. 
Suppose further that the previous polynomials belong to the same multihomogeneous component. For any $P$ let $\alpha(y_i):=\deg_{y_i}P$ and $\beta(z_i):=\deg_{z_i}P$. Then we consider the following substitution, to say, $\varphi$ such that
\[y_1\mapsto (e_{k+1}x_1)+e_{k+2}e_{k+3}+\cdots+e_{k+2(\alpha(y_1)-1)}e_{k+2\alpha(y_i)-1},\]
\[y_2\mapsto (e_{k+2\alpha(y_1)}x_2)+e_{k+2\alpha(y_1)+1}e_{k+2\alpha(y_i)+2}+\cdots+e_{k+2\alpha(y_1)+2\alpha(y_2)-3}e_{k+2\alpha(y_1)+2\alpha(y_2)-2},\]
\[\vdots\]
\[z_1\mapsto (e_1y_1)+e_2e_v+e_3e_{v+1}+\cdots+e_{\beta(z_1)}e_{v+\beta(z_1)-2},\]
\[z_2\mapsto (e_{\beta(z_1)+1}y_2)+e_{\beta(z_1)+2}e_{v+\beta(z_1)-1}+\cdots+e_{\beta(z_1)+\beta(z_2)}e_{v+\beta(z_1)+\beta(z_2)-3},\]
\[\vdots,\]
where $v$ is strictly greater than the maximum index of the $e_j$'s used for the $\varphi(y_i)$'s and the $x_i$'s and the $y_j$'s are basis elements of degree 0 of odd length with pairwise different supports. Unless explicitely written, the term inside the parenthesis has to be considered without the $x_i$ term.
Notice that we may always choose such a $v$ because we have an infinite number of $e_i$'s of $\Z_2$-degree 0. Moreover $\varphi$ is an admissible substitution because $\sum_i\beta(z_i)<k$. By Lemma \ref{useful1}, $\varphi(P)\neq0$ for each $P$. We consider the substitution $\psi$ obtained by $\varphi$ taking into account each parenthesis with the $x_i$ term. We use $\psi$ in order to annihilate the polynomials with a non-trivial commutator part. Then $\psi$ is non-zero only for one polynomial $P$ and this forces $\gamma(P)=0$.
We consider now $\psi$ obtained by $\varphi$ taking into account the $x_i$ term when the related homogeneous variable appears in the commutator part of one, precedently choosed, of the $P$'s. Then $\psi$ is non-zero only for $P$ and this forces $\gamma(P)=0$. We continue this process for each $P$. Hence we can prove that every subset of multihomogeneous polynomials of $S_1(n)$ is linearly independent. Suppose we have a linear combination $\sum_P\gamma(P)P\equiv0$ of elements of $S_1(n)$. 
Suppose further that the previous polynomials are in the same set of variables (because conversely the claim is trivial), then we consider again the substitution $\varphi$, where we consider $\alpha(y_i)$ to be the minimum degree of the $y_i$ appearing in the $P$'s such that $\deg_{y_j}P\leq\alpha(y_j)$ for each $j<i$ and $\beta(z_i)$ to be the minimum degree of the $z_i$ appearing in the $P$'s such that $\deg_{z_j}P\leq\beta(z_j)$ for each $j<i$ in order to annihilate all the polynomials having degree greater than $\alpha(y_i)$ in $y_i$ or greater than $\beta(z_i)$ in $z_i$. Then we reduce to the multihomogeneous case, hence $\gamma(P)=0$ for each $P$ and we are done.

In the case $F$ is of characteristic $p>0$ we observe that the only problems arise when dealing with the part of even variables. In fact in the case we have a variable $y$ of total degree greater than $p$, the substitution $\varphi$ vanishes on $y$ if $y$ has degree greater than $p$ and we cannot apply the machinery as above. In light of Lemma \ref{regevcons3} we avoid that problem using another substitution $\varphi_1$ such that $\varphi_1(y_i)\mapsto\varphi(y_i)+1$. Now the case in which all the polynomials are homogeneous follows by the first part of the proof. Suppose we have a finite sum $\sum_P\gamma(P)P\equiv0$ of polynomials that do not belong to the same multihomogeneous component but they are in the same set of variables. We work with polynomials without commutator part because the other case is analogous. Let $y_1$ be a variable appearing in all of the $P$'s, then we consider the following substitution \[\text{\rm$\psi_1(y_i)\mapsto1$ if $i\neq1$}\]\[\psi_1(y_1)\mapsto(1+a),\]where $a$ is an even element of $E$ as in Remark \ref{possiblechoose}. By Corollary \ref{regevcons2} we have that it suffices to consider only the case in which all of the exponents of $y_1$ are equal modulo $p$. Then the result follows by Lemma \ref{regevcons3}.
\endproof

We recall the following combinatorial tool. See the book of Stanley \cite{sta1} for more details.

\begin{proposition}\label{comb}
The number of commutative polynomials of degree $n$ in $k$ variables such that each variables has degree strictly less than $j$ is given by \[\kappa(n,j,k)=\sum_{r+sj=n}(-1)^s{k+r-1\choose r}{k\choose s}.\]
\end{proposition}

\begin{theorem}\label{dimension1}  Let $k\in\N$ and $p>2$ be a prime. The $\Z_2$-graded GK dimension of $E_{k^*}$ in $m$ variables over a field of characteristic 0 or $p$ is \[\gk_m^{\Z_2}(E_{k^*})=m.\]
\end{theorem}
\proof
We consider firstly the case of characteristic 0 or $p>k$. In light of Lemma \ref{independence}, each polynomial of $F_k^{\Z_2}(E_{k^*})$ can be written as a linear combination of $S_1$-type or $S_2$-type polynomials. Then in order to compute the $\Z_2$-graded GK dimension of $E_{k^*}$ in $m$ variables we have to count how many linearly independent elements lie in $S_1$ and $S_2$. We will count the number of polynomials of $S_1(t)$, then we argue analogously for $S_2(t)$. We observe that for any product $p_l$ of $l$ commutators of even (odd) variables we have ${t-2l+m-1\choose m-1}$ monomials $w$ such that $wp_l\in S_1(t)$. Due to the fact that we have ${m\choose 2l}$ choises for a product of $l$ commutators and that the number of commutators is less than $\left[\frac{m}{2}\right]$, we have that the total number of linearly independent polynomials of $S_1(t)$ is \tiny\[a_{k^*}(t)=\sum_{t_1+t_2=t,t_2\leq k}\sum_{l_1=0}^{\left[\frac{m}{2}\right]}\left({m\choose 2l_1}{t_1-2l_1+m-1\choose m-1}\right)\cdot\left(\sum_{l_2=0}^{\left[\frac{m}{2}\right]}{m\choose 2l_2}{t_2-2l_2+m-1\choose m-1}\right),\]\normalsize that is a polynomial in $t$ of degree $m-1$. Analogously the total number of linearly independent polynomials of $S_2(t)$ is \tiny\[b_{k^*}(t)=\sum_{t_1+t_2=t,t_2\leq k}m^2\left(\sum_{l_1=0}^{\left[\frac{m}{2}\right]}{m-1\choose 2l_1}{t_1-2l_1+m-2\choose m-1}\right)\cdot\left(\sum_{l_2=0}^{\left[\frac{m}{2}\right]}{m-1\choose 2l_2}{t_2-2l_2+m-2\choose m-1}\right),\]\normalsize where $m^2$ is the number of different $[z,y]$, that is a polynomial in $t$ of degree $m-1$. It turns out that the growth function of $E_{k^*}$ is a polynomial of degree $m$ and we are done. 

In the case the field is of characteristic $p\leq k$, we use Proposition \ref{comb} and we obtain that the total number of linearly independent polynomials of $S_1(t)$ is \tiny\[c_{k^*}(t)=\sum_{t_1+t_2=t,t_2\leq k}\sum_{l_1=0}^{\left[\frac{m}{2}\right]}\left({m\choose 2l_1}{t_1-2l_1+m-1\choose m-1}\right)\cdot\left(\sum_{l_2=0}^{\left[\frac{m}{2}\right]}{m\choose 2l_2}\kappa(t_2-2l_2,p,m)\right).\]\normalsize Analogously the total number of linearly independent polynomials of $S_2(t)$ is \tiny\[d_{k^*}(t)=\sum_{t_1+t_2=t,t_2\leq k}m^2\left(\sum_{l_1=0}^{\left[\frac{m}{2}\right]}{m-1\choose 2l_1}{t_1-2l_1+m-2\choose m-1}\right)\cdot\left(\sum_{l_2=0}^{\left[\frac{m}{2}\right]}{m-1\choose 2l_2}\kappa(t_2-2l_2-1,p,m)\right).\]\normalsize Then we have just to observe that the odd part is limited by a polynomial of degree $k$ and it does not contribute to the total GK dimension and we are done.
\endproof

\begin{corollary}
The $\Z_2$-graded Hilbert series of $E_{k^*}$ is \[H(E_{k^*};T)=\sum_{t\geq0}(a_{k^*}(t)+b_{k^*}(t))T^t,\] if $F$ is a field of characteristic 0 or $p>k$. 
The $\Z_2$-graded Hilbert series of $E_{k^*}$ is \[H(E_{k^*};T)=\sum_{t\geq0}(c_{k^*}(t)+d_{k^*}(t))T^t,\] if $F$ is a field of characteristic $p\leq k$.
\end{corollary}

\section{The graded GK dimension of $E_\infty$}

Let $p>2$ be a prime and $E_{\infty}$ as defined above over a field of characteristic 0.
In light of Thorem 10 of \cite{did1}, $E_\infty$ satisfies the identity $[x_1,x_2,x_3]$, then each polynomial of the relatively-free $\Z_2$-graded algebra of $E_\infty$ can be written in the following form: \[y_{i_1}^{\alpha_{i_1}}\cdots y_{i_l}^{\alpha_{i_l}}z_{j_1}^{\beta_{j_1}}\cdots z_{j_r}^{\beta_{j_r}}[z_{s_1},z_{s_2}]\cdots[z_{s_p},z_{s_{p+1}}][y_{t_1},y_{t_2}]\cdots[y_{t_q},y_{t_{q+1}}]\]
\[y_{i_1}^{\alpha_{i_1}}\cdots y_{i_l}^{\alpha_{i_l}}z_{j_1}^{\beta_{j_1}}\cdots z_{j_r}^{\beta_{j_r}}[z_{s_1},z_{s_2}]\cdots[z_{s_p},z_{s_{p+1}}][z_{s_{p+2}},y_{t_1}][y_{t_2},y_{t_3}]\cdots[y_{t_q},y_{t_{q+1}}]\] without any restriction on the $z$'s. In the case $F$ is a field of characteristic $p$ we may suppose each $\beta_{j_l}$ to be strictly less than $p$. Following word by word the proof of Lemma \ref{independence} we have the following.

\begin{lemma}\label{independence2}
The polynomials of the following type: \[y_{i_1}^{\alpha_{i_1}}\cdots y_{i_l}^{\alpha_{i_l}}z_{j_1}^{\beta_{j_1}}\cdots z_{j_m}^{\beta_{j_m}}[z_{s_1},z_{s_2}]\cdots[z_{s_p},z_{s_{p+1}}][y_{t_1},y_{t_2}]\cdots[y_{t_q},y_{t_{q+1}}]\]
\[y_{i_1}^{\alpha_{i_1}}\cdots y_{i_l}^{\alpha_{i_l}}z_{j_1}^{\beta_{j_1}}\cdots z_{j_r}^{\beta_{j_r}}[z_{s_1},z_{s_2}]\cdots[z_{s_p},z_{s_{p+1}}][z_{s_{p+2}},y_{t_1}][y_{t_2},y_{t_3}]\cdots[y_{t_q},y_{t_{q+1}}]\]
are linearly independent modulo $T_{\Z_2}(E_{\infty})$.
\end{lemma}

We observe that the substitution done for $E_{k^*}$ is still valid in this case due to the fact that $E_\infty$ has an infinite number of $e_i$'s of $\Z_2$-degree 0.

\begin{theorem}\label{dimension2}Let $k\in\N$ and $p>2$ be a prime. Then the $\Z_2$-graded GK dimension of $E_{\infty}$ in $m$ variables is \[\text{$\gk_m^{\Z_2}(E_{\infty})=2m$ if the characteristic of $F$ is 0},\]
\[\text{$\gk_m^{\Z_2}(E_{\infty})=m$ if the characteristic of $F$ is $p$}.\]
\end{theorem}
\proof
We consider firstly the case of characteristic 0. In light of Lemma \ref{independence2} we have that the polynomials \[y_{i_1}^{\alpha_{i_1}}\cdots y_{i_l}^{\alpha_{i_l}}z_{j_1}^{\beta_{j_1}}\cdots z_{j_r}^{\beta_{j_r}}[z_{s_1},z_{s_2}]\cdots[z_{s_p},z_{s_{p+1}}][y_{t_1},y_{t_2}]\cdots[y_{t_q},y_{t_{q+1}}],\]
\[y_{i_1}^{\alpha_{i_1}}\cdots y_{i_l}^{\alpha_{i_l}}z_{j_1}^{\beta_{j_1}}\cdots z_{j_r}^{\beta_{j_r}}[z_{s_1},z_{s_2}]\cdots[z_{s_p},z_{s_{p+1}}][z_{s_{p+2}},y_{t_1}][y_{t_2},y_{t_3}]\cdots[y_{t_q},y_{t_{q+1}}]\] are linearly independent. We compute their number using the same argument of Theorem \ref{dimension1}. We have that the total number of linearly independent polynomials of $S_1(t)$ is \tiny\[a_{\infty}(t)=\sum_{t_1+t_2=t}\left(\sum_{l_1=0}^{\left[\frac{m}{2}\right]}{m\choose 2l_1}{t_1-2l_1+m-1\choose m-1}\right)\cdot\left(\sum_{l_2=0}^{\left[\frac{m}{2}\right]}{m\choose 2l_2}{t_2-2l_2+m-1\choose m-1}\right),\]\normalsize that is a polynomial in $t$ of degree $2m-1$. Analogously the total number of linearly independent polynomials of $S_2(t)$ is \tiny\[b_{\infty}(t)=\sum_{t_1+t_2=t}m^2\left(\sum_{l_1=0}^{\left[\frac{m}{2}\right]}{m-1\choose 2l_1}{t_1-2l_1+m-2\choose m-1}\right)\cdot\left(\sum_{l_2=0}^{\left[\frac{m}{2}\right]}{m-1\choose 2l_2}{t_2-2l_2+m-2\choose m-1}\right),\]\normalsize that is a polynomial in $t$ of degree $2m-1$. It turns out that the growth function of $E_{\infty}$ is a polynomial of degree $2m$ and we are done. In the case $F$ has characteristic $p$, we use Proposition \ref{comb}, then the total number of linearly independent polynomials of $S_1(t)$ is \tiny\[c_{\infty}(t)=\sum_{t_1+t_2=t}\left(\sum_{l_1=0}^{\left[\frac{m}{2}\right]}{m\choose 2l_1}{t_1-2l_1+m-1\choose m-1}\right)\cdot\left(\sum_{l_2=0}^{\left[\frac{m}{2}\right]}{m\choose 2l_2}\kappa(t_2-2l_2,p,m)\right).\]\normalsize Analogously the total number of linearly independent polynomials of $S_2(t)$ is \tiny\[d_{\infty}(t)=\sum_{t_1+t_2=t}m^2\left(\sum_{l_1=0}^{\left[\frac{m}{2}\right]}{m-1\choose 2l_1}{t_1-2l_1+m-2\choose m-1}\right)\cdot\left(\sum_{l_2=0}^{\left[\frac{m}{2}\right]}{m-1\choose 2l_2}\kappa(t_2-2l_2-1,p,m)\right).\]\normalsize We observe that the $z$ part is limited by a polynomial of degree $m(p-1)+\left[\frac{m}{2}\right]$ and we are done.
\endproof

\begin{remark}
We observe that if $F$ is a field of characteristic 0, then the $T_{\Z_2}$ ideal of $E_\infty$ ``coincides'' with the $T$-ideal of $E$. Nevertheless the $\Z_2$-graded GK dimension of $E_\infty$ is the GK dimension of $F_m^{\Z_2}(E_\infty)$ that is merely the GK dimension of a relatively-free algebra of $2m$ variables ``modulo'' the identity $[u_1,u_2,u_3]$. It turns out it must be equal to the GK dimension of $E$ in $2m$ variables as we may verify by \cite{cen2,cen4}.
\end{remark}

\begin{corollary}
The $\Z_2$-graded Hilbert series of $E_{\infty}$ is \[H(E_{\infty};T)=\sum_{t\geq0}(a_{\infty}(t)+b_{\infty}(t))T^t,\] if $F$ is a field of characteristic 0. 
The $\Z_2$-graded Hilbert series of $E_{\infty}$ is \[H(E_{\infty};T)=\sum_{t\geq0}(c_{\infty}(t)+d_{\infty}(t))T^t,\] if $F$ is a field of characteristic $p$.
\end{corollary}

\section{The graded GK dimension of $E_k$}

Let us consider $k\in\N$ and $p>2$ be a prime such that $p>k$ and $E_{k}$ as defined above over a field of characteristic 0 or $p>k$. Suppose $k$ is even, then we argue analogously for $k$ odd. By Theorem 38 of \cite{did1} and by Theorem \ref{t3} we have \[[y_1,y_2][y_3,y_4]\cdots[y_{k-1},y_k][y_{k+1},u_{k+2}],\] where $u_{k+2}$ is any homogenous variable, is a $\Z_2$-graded identities for $E_k$. Working inside $F^{\Z_2}_m(E_k)$ we note that we may consider only polynomials with at most $s:=\left[\frac{k}{2}\right]$ commutators in the $y$'s. By Theorem \ref{t3} we have \[g_{k-l+2}(z_1,\ldots,z_{k-l+2})[y_1,y_2]\cdots[y_{l-1},y_l],\] where $l\leq k$, is a $\Z_2$-graded identity for $E_k$. If $l=k$, then $g_2(z_1,z_2)=z_1z_2-\frac{1}{2}[z_1,z_2]$. Hence \[z_1z_2[y_1,y_2]\cdots[y_{k-1},y_k]\equiv\frac{1}{2}[z_1,z_2][y_1,y_2]\cdots[y_{k-1},y_k].\]In general, each polynomial may be written modulo $T_{\Z_2}(E_k)$ as a linear combination of the following polynomials 
\[y_{i_1}^{\alpha_{i_1}}\cdots y_{i_l}^{\alpha_{i_l}}z_{j_1}^{\beta_{j_1}}\cdots z_{j_r}^{\beta_{j_r}}[z_{s_1},z_{s_2}]\cdots[z_{s_p},z_{s_{p+1}}][y_{t_1},y_{t_2}]\cdots[y_{t_q},y_{t_{q+1}}]\]
\[y_{i_1}^{\alpha_{i_1}}\cdots y_{i_l}^{\alpha_{i_l}}z_{j_1}^{\beta_{j_1}}\cdots z_{j_r}^{\beta_{j_r}}[z_{s_1},z_{s_2}]\cdots[z_{s_p},z_{s_{p+1}}][z_{s_{p+2}},y_{t_1}][y_{t_2},y_{t_3}]\cdots[y_{t_q},y_{t_{q+1}}],\] where the polynomial part in the $z$'s is a polynomial of total degree strictly less than $k-(q+1)+1$.

\begin{lemma}\label{useit}
The polynomials of the following type: 
\[y_{i_1}^{\alpha_{i_1}}\cdots y_{i_l}^{\alpha_{i_l}}z_{j_1}^{\beta_{j_1}}\cdots z_{j_r}^{\beta_{j_r}}[z_{s_1},z_{s_2}]\cdots[z_{s_p},z_{s_{p+1}}][y_{t_1},y_{t_2}]\cdots[y_{t_q},y_{t_{q+1}}]\]
\[y_{i_1}^{\alpha_{i_1}}\cdots y_{i_l}^{\alpha_{i_l}}z_{j_1}^{\beta_{j_1}}\cdots z_{j_r}^{\beta_{j_r}}[z_{s_1},z_{s_2}]\cdots[z_{s_p},z_{s_{p+1}}][z_{s_{p+2}},y_{t_1}][y_{t_2},y_{t_3}]\cdots[y_{t_q},y_{t_{q+1}}],\] 
 where $q+1$ is fixed, are linearly independent modulo $T_{\Z_2}(E_{k})$.
\end{lemma}

\proof 
We will argue only for the $S_1$-type polynomials because in the other case the proof is analogous. Let us consider the set $S_1(n)$ of polynomials of $S_1$ \[y_{i_1}^{\alpha_{i_1}}\cdots y_{i_l}^{\alpha_{i_l}}z_{j_1}^{\beta_{j_1}}\cdots z_{j_r}^{\beta_{j_r}}[z_{s_1},z_{s_2}]\cdots[z_{s_{t-1}},z_{s_{t}}][y_{t_1},y_{t_2}]\cdots[y_{t_{q-1}},y_{t_{q}}]\] of total degree $n$. Suppose we have a linear combination $\sum_P\gamma(p)P\equiv0$ of elements of $S_1(n)$. 
Suppose further that the previous polynomials belong to the same multihomogeneous component. Let $\sum_P\gamma(P)P$ the linear combination of polynomials of $S_1(n)$ belonging to the same multihomogeneous component and we observe that the part in the $y_i$'s is completely determined by the variables appearing in the commutators. For any $P$ let $\alpha(y_i):=\deg_{z_i}P$ and $\beta(z_i):=\deg_{z_i}P$. Then we consider the following substitution, to say, $\varphi$:
\small\[y_1\mapsto (e_1)+e_{k+1}e_{k+2}+\cdots+e_{k+2(\alpha(y_1)-1)-1}e_{k+2(\alpha(y_1)-1)},\]
\small\[y_2\mapsto (e_2)+e_{k+2(\alpha(y_1)-1)+1}e_{k+2(\alpha(y_1)-1)+2}+\cdots+e_{k+2(\alpha(y_1)-1)+2(\alpha(y_2)-1)-1}e_{k+2(\alpha(y_1)-1)+2(\alpha(y_2)-1)},\]
\[\vdots\]
\small\[y_q\mapsto (e_q)+e_{k+\sum_{i=1}^{q-1} 2(\alpha(y_i)-1)+1}e_{k+\sum 2(\alpha(y_i)-1)+2}+\cdots+e_{k+\sum_{i=1}^{q-1} 2(\alpha(y_i)-1)-1}e_{k+\sum_{i=1}^{q-1} 2(\alpha(y_i)-1)},\]
\small\[y_{q+1}\mapsto e_{v+1}e_{v+2}+\cdots+e_{v+2(\alpha(y_{q+1})-1)-1}e_{v+2(\alpha(y_{q+1})-1)},\]
\[\vdots,\]
\small\[y_m\mapsto e_{v+\sum_{i=q+1}^{m-1} 2(\alpha(y_i)-1)+1}e_{v+\sum 2(\alpha(y_i)-1)+2}+\cdots+e_{v+\sum_{i=1}^{q-1} 2(\alpha(y_i)-1)-1}e_{v+\sum_{i=1}^{q-1} 2(\alpha(y_i)-1)},\]
\small\[z_1\mapsto (e_{u+1})+e_{q+1}e_{u+2}+\cdots+e_{q+\beta(z_1)-1}e_{u+\beta(z_1)},\]
\small\[z_2\mapsto (e_{u+\beta(z_1)+1})+e_{q+\beta(z_1)}e_{u+\beta(z_1)+2}+\cdots+e_{q+\beta(z_1)+\beta(z_2)-2}e_{u+\beta(z_1)+\beta(z_2)},\]
\[\vdots,\]
\small\[z_t\mapsto (e_{u+\sum_{i=1}^{t-1}\beta(z_i)+1})+e_{q+\sum_{i=1}^{t-1}\beta(z_i)-t+2}e_{u+\sum_{i=1}^{t-1}\beta(z_i)+2}\cdots+e_{q+\sum_{i=1}^{t}\beta(z_i)-t}e_{u+\sum_{i=1}^{t}\beta(z_i)},\]
\small\[z_{t+1}\mapsto e_{c+1},\]
\[\vdots,\]
\small\[z_m\mapsto e_{c+m-t},\]\normalsize
where $v=k+\sum_{i=1}^{q-1} 2(\alpha(y_i)-1)$, $u=v+\sum_{i=1}^{q-1} 2(\alpha(y_i)-1)$ and $c=u+\sum_{i=1}^{t}\beta(z_i)$. Notice that we may always choose $u,v,c$ because we have an infinite number of $e_i$'s of $\Z_2$-degree 1. Moreover $\varphi$ is an admissible substitution because $\sum_i\beta(z_i)<k-q+1$. We observe that by Lemma \ref{useful1}, $\varphi(P)\neq0$ for each $P$. In fact if there exists a variable $z$ in the polynomial part appearing in the commutator part too, we consider the substitution $\varphi(z)$ including the part in the parenthesis. Now the proof follows using the same stuffs as in Lemma \ref{independence}.
\endproof

\begin{theorem}Let $k\in\N$ and $p>2$ be a prime. Then the $\Z_2$-graded GK dimension of $E_{k}$ in $m$ variables is \[\text{$\gk_m^{\Z_2}(E_{k})=m$}.\]
\end{theorem}
\proof
Let $F$ be a field of characteristic 0 or $p>k$. We observe that the number of commutators in the $y$ is limited by $\left[\frac{k}{2}\right]$. We fix now $l\leq k$, $l$ even, then by Lemma \ref{useit} the polynomials \[y_{i_1}^{\alpha_{i_1}}\cdots y_{i_l}^{\alpha_{i_l}}z_{j_1}^{\beta_{j_1}}\cdots z_{j_m}^{\beta_{j_m}}[z_{s_1},z_{s_2}]\cdots[z_{s_p},z_{s_{p+1}}][y_{t_1},y_{t_2}]\cdots[y_{t_{l-1}},y_{t_{l}}]\]  are linearly independent and their number is \tiny\[a_{k}(t)=\sum_{{\scriptsize\begin{array}{cc}
	t_1+t_2=t\\
	\end{array}}}\sum_{{\scriptsize\begin{array}{cc}
	l\leq k\\
	\text{\rm $l$ even}
\end{array}}}\sum_{{\scriptsize\begin{array}{cc}
	s=0\\
	t_2-2s\leq \beta\end{array}}}^{\left[\frac{m}{2}\right]}{m\choose l}{t_1-l+m-1\choose m-1}{m\choose 2s}{t_2-2s+m-1\choose m-1},\]\normalsize where $\beta=k-l$. We argue analogously in the case $l$ is odd, then the total number of linearly independent polynomials in this case is \tiny\[b_{k}(t)=\sum_{{\scriptsize\begin{array}{cc}
	t_1+t_2=t\\
	\end{array}}}\sum_{{\scriptsize
\begin{array}{cc}
	l\leq k\\
	\text{\rm $l$ odd}
\end{array}}
}\sum_{{\scriptsize\begin{array}{cc}
	s=0\\
	t_2-2s\leq \beta\end{array}}}^{\left[\frac{m}{2}\right]}{m\choose l-1}m(m-l+1){t_1-l+m-1\choose m-1}{m-1\choose 2s}{t_2-2s+m-1\choose m-1}.\]\normalsize
Due to the fact that such polynomials are essentially polynomials in the $y$'s with a limited-by-polynomial part in the $z$'s, we are just allowed to take into account the monomials in the $y$'s. Hence $\gk_m^{\Z_2}(E_{k})=m$ and we are done. If the characteristic of $F$ is $p\leq k$, then the total number of the linear independent polynomials is \tiny\[c_{k}(t)=\sum_{{\scriptsize\begin{array}{cc}
	t_1+t_2=t\\
	\end{array}}}\sum_{{\scriptsize\begin{array}{cc}
	l\leq k\\
	\text{\rm $l$ even}
\end{array}}}\sum_{{\scriptsize\begin{array}{cc}
	s=0\\
	t_2-2s\leq \beta\end{array}}}^{\left[\frac{m}{2}\right]}{m\choose l}{t_1-l+m-1\choose m-1}{m\choose 2s}\kappa(t_2-2s,p,m-1),\]\normalsize where $\beta=k-l$. We argue analogously in the case $l$ is odd, then the total number of linearly independent polynomials is \tiny\[d_{k}(t)=\sum_{{\scriptsize\begin{array}{cc}
	t_1+t_2=t\\
	\end{array}}}\sum_{{\scriptsize
\begin{array}{cc}
	l\leq k\\
	\text{\rm $l$ odd}
\end{array}}
}\sum_{{\scriptsize\begin{array}{cc}
	s=0\\
	t_2-2s\leq \beta\end{array}}}^{\left[\frac{m}{2}\right]}{m\choose l-1}m(m-l+1){t_1-l+m-1\choose m-1}{m-1\choose 2s}\kappa(t_2-2s,p,m-1),\]\normalsize then we argue as in the previous case.
\endproof

\begin{corollary}
The $\Z_2$-graded Hilbert series of $E_{k}$ is \[H(E_{k^*};T)=\sum_{t\geq0}(a_{k}(t)+b_{k}(t))T^t,\] if $F$ is a field of characteristic 0 or $p>k$. 
The $\Z_2$-graded Hilbert series of $E_{k}$ is \[H(E_{k^*};T)=\sum_{t\geq0}(c_{k}(t)+d_{k}(t))T^t,\] if $F$ is a field of characteristic $p\leq k$.
\end{corollary}

\section{Conclusions}

We recall the definition of \textit{exponent} of a PI-algebra. Let $A$ be a PI-algebra and let $n\in\N$. If $V_n$ is the $F$-space of multilinear polynomials in $n$ variables we may study the factor space $V_n(A):=V_n/(V_n\cap T(A))$. We observe that $V_n(A)$ is a vector space over $F$, then we may take into an account its dimension $c_n(A)$, called $n$\textit{-th codimension} of $A$. In \cite{reg1} Regev proved that the sequence of codimensions of a PI-algebra is exponentially bounded. In \cite{giz1}, \cite{giz2} Giambruno and Zaicev proved that:

\begin{theorem}\label{teogiambrunozaicev} If $A$ is any PI-algebra, then there exists the limit
\[\text{\rm exp}(A)=\lim_{n\rightarrow\infty}\sqrt[n]{c_n(A)}\] and it is a non-negative integer called the PI-exponent of $A$.\end{theorem}

We consider now a $G$-graded algebra and let us denote by $1_G$ the identity of the group $G$. Some of the relations between the identities of $A^{1_G}$ and $A$ have been studied in a number of papers. Speaking in terms of exponent, it is rather easy to see that $\ep(A^{1_G})\leq\ep(A)$ since $A^{1_G}$ is a subalgebra of $A$. On the other hand, there are a lot of examples of graded algebras satisfying the inequality \begin{equation}\label{aljadeff}\ep(A)\leq\left|G\right|^2\ep(A^{1_G}).\end{equation} In \cite{baz1}, Bahturin and Zaicev conjectured that the inequality (\ref{aljadeff}) holds for any graded algebra with a non-trivial polynomial identity. In the case $G$ is a fnite abelian group, the pevious conjecture has been already solved positively by Aljadeff, Giambruno and La Mattina (see \cite{agl1}).

In the spirit of \ref{aljadeff} and in the light of the result of \cite{agl1}, Aljadeff and the author conjectured that:

\begin{conjecture}\label{alc}
Let $G$ be a finite abelian group and $A$ be a $G$-graded PI-algebra, then for any $k\geq1$ \[\gk_k^G(A)\leq|G|^2\gk_k(A^{1_G}).\]
\end{conjecture}

We observe that in each of the cases we have investigated, the even part $E_d^0$ grows as a polynomial algebra in $m$ variables, hence $\gk_m(E_d^0)=m$. Moreover notice that the maximal $\Z_2$-graded GK dimension is $2m$ and is obtained in the case $E$ is graded by $\deg_\infty$. Hence $\gk_m^{\Z_2}(E)\leq 2m$, then a fortiori it satisfies Conjecture \ref{alc}. Clearly in this case the term $|G|^2$ is too much but in a forthcoming paper Aljadeff and the author show that $|G|^2$ is needed when working with matrices graded by finite abelian groups of order $n$.

\end{document}